\documentclass{article}
\usepackage{amsmath}

\newtheorem{theorem}{Theorem}

\newtheorem{corollary}[theorem]{Corollary}

\newtheorem{lemma}[theorem]{Lemma}

\newtheorem{remark}[theorem]{Remark}

\newenvironment{proof}[1][Proof]{\textbf{#1.} }{\ \rule{0.5em}{0.5em}}

\begin{document}

\title{Peano-like bounds for some Newton-Cotes Formulae}
\author{Nenad Ujevi\'{c} \\
Department of Mathematics\\
University of Split\\
Teslina 12/III, 21000 Split\\
CROATIA}
\maketitle

\begin{abstract}
An error analysis for some Newton-Cotes quadrature formulae is presented.
Peano-like error bounds are obtained. They are generally, but not always,
better than the usual Peano bounds.
\end{abstract}



\section{Introduction}

In this paper we present an error analysis for some Newton-Cotes quadrature
formulae. We consider Simpson's rule, 3/8 Simpson rule and Boole's rule. A
similar error analysis for Simpson's rule have been investigated more
recently (\cite{C1}, \cite{DAC1}, \cite{DCR1}, \cite{DPW1}, \cite{PPUV1})
with the view of obtaining bounds on the quadrature rule in terms of a
variety of norms involving, at most, the first derivative. It is well known
that if the mapping $f$ is neither four times differentiable nor is the
fourth derivative $f^{(4)}$ bounded, then we cannot apply the classical
Simpson's quadrature formula, which, actually is one of the most used
quadrature formulas in practical applications. Thus, the above mentioned
analysis is important as well as the analysis presented here.

The current work brings results for the above mentioned Newton-Cotes
quadrature rules giving explicit error bounds and using results from the
modern theory of inequalities. The used inequalities are known in the
literature as inequalities of Ostrowski-Gr\"{u}ss type. The error bounds are
expressed in terms of second derivatives. As we have already mentioned for
Simpson's rule, the general approach used in the past involves the
assumption of bounded derivatives of degree higher than two. We also mention
that the obtained results can be derived using the Peano kernel theorem. In
any case, these bounds are generally, but not always, better that the usual
Peano error bounds (see Remarks \ref{R1}, \ref{R2} and \ref{R3}).

Here we do not consider composite quadrature rules since they can be formed
in the usual way. However, the analysis presented here allows the
determination of the partition required that would assure the accuracy the
result would be within a prescribed error tolerance.

In Section 2 we established some auxiliary results. We use these results in
further sections. In Section 3 we consider Simpson's rule. In Section 4 we
consider 3/8 Simpson rule and in Section 5 we consider Boole's rule.

\section{Preliminary results}

\begin{lemma}
\label{L1}Let $I\subset R$ be an open interval and $a,b\in I,$ $a<b.$ Let $%
f:I\rightarrow R$ be a twice differentiable function and $\ $let $x\in \left[
a,b\right] $ be a fixed element. Then we have%
\begin{eqnarray}
&&f(x)(b-a)-(x-\frac{a+b}{2})\left[ f(b)-f(a)\right] -\int%
\limits_{a}^{b}f(t)dt  \label{j1} \\
&=&\frac{1}{b-a}\int\limits_{a}^{b}\int\limits_{a}^{b}p(x,t)p(t,s)f^{\prime
\prime }(s)dsdt,  \notag
\end{eqnarray}%
where%
\begin{equation}
p(x,t)=\left\{ 
\begin{array}{c}
t-a,\quad t\in \left[ a,x\right] \\ 
t-b,\quad t\in \left( x,b\right]%
\end{array}%
\right. .  \label{j2}
\end{equation}
\end{lemma}

\begin{proof}
Integrating by parts, we have%
\begin{equation*}
\int\limits_{a}^{b}p(x,t)f^{\prime
}(t)dt=f(x)(b-a)-\int\limits_{a}^{b}f(t)dt,
\end{equation*}%
i. e.%
\begin{equation*}
f(x)=\frac{1}{b-a}\int\limits_{a}^{b}f(t)dt+\frac{1}{b-a}\int%
\limits_{a}^{b}p(x,t)f^{\prime }(t)dt.
\end{equation*}%
If we substitute $f\rightarrow f^{\prime }$ in the above relation, then we
get%
\begin{equation*}
f^{\prime }(t)=\frac{1}{b-a}\left[ f(b)-f(a)\right] +\frac{1}{b-a}%
\int\limits_{a}^{b}p(t,s)f^{\prime \prime }(s)ds.
\end{equation*}%
Thus,%
\begin{eqnarray*}
&&\int\limits_{a}^{b}p(x,t)f^{\prime }(t)dt \\
&=&\int\limits_{a}^{b}p(x,t)\left[ \frac{1}{b-a}\left[ f(b)-f(a)\right] +%
\frac{1}{b-a}\int\limits_{a}^{b}p(t,s)f^{\prime \prime }(s)ds\right] dt.
\end{eqnarray*}%
We also have%
\begin{equation*}
\int\limits_{a}^{b}p(x,t)dt=(b-a)(x-\frac{a+b}{2}).
\end{equation*}%
From the above relations it follows%
\begin{eqnarray*}
&&\frac{1}{b-a}\int\limits_{a}^{b}\int\limits_{a}^{b}p(x,t)p(t,s)f^{\prime
\prime }(s)dsdt \\
&=&f(x)(b-a)-(x-\frac{a+b}{2})\left[ f(b)-f(a)\right] -\int%
\limits_{a}^{b}f(t)dt.
\end{eqnarray*}
\end{proof}

\begin{lemma}
\label{L2}Let $p(x,t)$ be defined by (\ref{j2}). Then we have%
\begin{eqnarray}
q(x,s) &=&\int\limits_{a}^{b}p(x,t)p(t,s)dt  \label{j3} \\
&=&\left\{ 
\begin{array}{c}
(b-a)(x-\frac{a+b}{2})(s-a)-\frac{b-a}{2}(s-a)^{2},\quad s\in \left[ a,x%
\right] \\ 
(b-a)(x-\frac{a+b}{2})(s-b)-\frac{b-a}{2}(s-b)^{2},\quad s\in \left( x,b%
\right]%
\end{array}%
\right. ,  \notag
\end{eqnarray}%
where $x\in \left[ a,b\right] .$
\end{lemma}

\begin{proof}
For $s\in \left[ a,x\right] $ we have%
\begin{eqnarray*}
q(x,s)
&=&\int\limits_{a}^{s}(t-a)(s-b)dt+\int\limits_{s}^{x}(t-a)(s-a)dt+\int%
\limits_{x}^{b}(t-b)(s-a)dt \\
&=&(s-b)\frac{(s-a)^{2}}{2}+(s-a)\frac{(x-a)^{2}-(s-a)^{2}}{2}-(s-a)\frac{%
(x-b)^{2}}{2} \\
&=&(b-a)(x-\frac{a+b}{2})(s-a)-\frac{b-a}{2}(s-a)^{2}.
\end{eqnarray*}%
For $s\in \left( x,b\right] $ we have%
\begin{eqnarray*}
q(x,s)
&=&\int\limits_{a}^{x}(t-a)(s-b)dt+\int\limits_{x}^{s}(t-b)(s-b)dt+\int%
\limits_{s}^{b}(t-b)(s-a)dt \\
&=&(s-a)\frac{(x-a)^{2}}{2}+(s-b)\frac{(x-b)^{2}+(s-b)^{2}}{2}-(s-a)\frac{%
(s-b)^{2}}{2} \\
&=&(b-a)(x-\frac{a+b}{2})(s-b)-\frac{b-a}{2}(s-b)^{2}.
\end{eqnarray*}%
From the above relations we see that (\ref{j3}) holds.
\end{proof}

\begin{corollary}
\label{C1}Let the assumptions of Lemma \ref{L2} be satisfied. Then we have%
\begin{eqnarray*}
q(a,s) &=&-\frac{1}{2}(b-a)^{2}(s-b)-\frac{b-a}{2}(s-b)^{2}, \\
q(b,s) &=&\frac{1}{2}(b-a)^{2}(s-a)-\frac{b-a}{2}(s-a)^{2}, \\
q(\frac{a+b}{2},s) &=&\left\{ 
\begin{array}{c}
-\frac{b-a}{2}(s-a)^{2},\quad s\in \left[ a,\frac{a+b}{2}\right] \\ 
-\frac{b-a}{2}(s-b)^{2},\quad s\in \left( \frac{a+b}{2},b\right]%
\end{array}%
\right. , \\
q(\frac{3a+b}{4},s) &=&\left\{ 
\begin{array}{c}
-\frac{1}{4}(b-a)^{2}(s-a)-\frac{b-a}{2}(s-a)^{2},\quad s\in \left[ a,\frac{%
3a+b}{4}\right] \\ 
-\frac{1}{4}(b-a)^{2}(s-b)-\frac{b-a}{2}(s-b)^{2},\quad s\in \left( \frac{%
3a+b}{4},b\right]%
\end{array}%
\right. , \\
q(\frac{a+3b}{4},s) &=&\left\{ 
\begin{array}{c}
\frac{1}{4}(b-a)^{2}(s-a)-\frac{b-a}{2}(s-a)^{2},\quad s\in \left[ a,\frac{%
a+3b}{4}\right] \\ 
\frac{1}{4}(b-a)^{2}(s-b)-\frac{b-a}{2}(s-b)^{2},\quad s\in \left( \frac{a+3b%
}{4},b\right]%
\end{array}%
\right. , \\
q(\frac{2a+b}{3},s) &=&\left\{ 
\begin{array}{c}
-\frac{1}{6}(b-a)^{2}(s-a)-\frac{b-a}{2}(s-a)^{2},\quad s\in \left[ a,\frac{%
2a+b}{3}\right] \\ 
-\frac{1}{6}(b-a)^{2}(s-b)-\frac{b-a}{2}(s-b)^{2},\quad s\in \left( \frac{%
2a+b}{3},b\right]%
\end{array}%
\right. , \\
q(\frac{a+2b}{3},s) &=&\left\{ 
\begin{array}{c}
\frac{1}{6}(b-a)^{2}(s-a)-\frac{b-a}{2}(s-a)^{2},\quad s\in \left[ a,\frac{%
a+2b}{3}\right] \\ 
\frac{1}{6}(b-a)^{2}(s-b)-\frac{b-a}{2}(s-b)^{2},\quad s\in \left( \frac{a+2b%
}{3},b\right]%
\end{array}%
\right. .
\end{eqnarray*}
\end{corollary}

\section{Simpson's rule}

\begin{theorem}
\label{T1}Let the assumptions of Lemma \ref{L1} hold and let $\gamma ,\Gamma 
$ be real numbers such that $\gamma \leq f^{\prime \prime }(t)\leq \Gamma ,$ 
$\forall t\in \left[ a,b\right] $.\ Then we have%
\begin{equation}
\left| \frac{b-a}{6}\left[ f(a)+4f(\frac{a+b}{2})+f(b)\right]
-\int\limits_{a}^{b}f(t)dt\right| \leq \frac{\Gamma -\gamma }{162}(b-a)^{3}.
\label{j7}
\end{equation}
\end{theorem}

\begin{proof}
For $x=a$ the left-hand side in (\ref{j1}) is equal to%
\begin{equation}
(b-a)\frac{f(b)+f(a)}{2}-\int\limits_{a}^{b}f(t)dt.  \label{j8}
\end{equation}%
For $x=b$ the left-hand side in (\ref{j1}) is equal to%
\begin{equation}
(b-a)\frac{f(b)+f(a)}{2}-\int\limits_{a}^{b}f(t)dt.  \label{j9}
\end{equation}%
For $x=\frac{a+b}{2}$ the left-hand side in (\ref{j1}) is equal to%
\begin{equation}
(b-a)f(\frac{a+b}{2})-\int\limits_{a}^{b}f(t)dt.  \label{j10}
\end{equation}%
If we now multiply (\ref{j10}) by 2 and add (\ref{j8}) or (\ref{j9}) then we
get%
\begin{equation}
\frac{b-a}{2}\left[ f(a)+4f(\frac{a+b}{2})+f(b)\right] -3\int%
\limits_{a}^{b}f(t)dt.  \label{j11a}
\end{equation}%
The corresponding right-hand side is%
\begin{eqnarray}
R(a,b) &=&\frac{1}{b-a}\int\limits_{a}^{b}\int\limits_{a}^{b}\left[ 2p(\frac{%
a+b}{2},t)+p(a,t)\right] p(t,s)f^{\prime \prime }(s)dtds  \label{j13} \\
&=&\frac{1}{b-a}\int\limits_{a}^{b}\left[ 2q(\frac{a+b}{2},s)+q(a,s)\right]
f^{\prime \prime }(s)ds.  \notag
\end{eqnarray}%
Let $K_{1}(s)=2q(\frac{a+b}{2},s)+q(a,s)$. Then we have%
\begin{equation*}
\int\limits_{a}^{b}K_{1}(s)ds=0
\end{equation*}%
such that%
\begin{equation*}
\frac{1}{b-a}\int\limits_{a}^{b}K_{1}(s)f^{\prime \prime }(s)ds=\frac{1}{b-a}%
\int\limits_{a}^{b}K_{1}(s)\left[ f^{\prime \prime }(s)-\frac{\Gamma +\gamma 
}{2}\right] ds
\end{equation*}%
and%
\begin{eqnarray*}
\left| R(a,b)\right| &\leq &\frac{1}{b-a}\underset{s\in \left[ a,b\right] }{%
\max }\left| f^{\prime \prime }(s)-\frac{\Gamma +\gamma }{2}\right|
\int\limits_{a}^{b}\left| K_{1}(s)\right| ds \\
&\leq &\frac{\Gamma -\gamma }{2(b-a)}\int\limits_{a}^{b}\left|
K_{1}(s)\right| ds,
\end{eqnarray*}%
since%
\begin{equation}
\underset{s\in \left[ a,b\right] }{\max }\left| f^{\prime \prime }(s)-\frac{%
\Gamma +\gamma }{2}\right| \leq \frac{\Gamma -\gamma }{2}.  \label{zvj}
\end{equation}

Hence,%
\begin{equation}
\left| R(a,b)\right| \leq \frac{\Gamma -\gamma }{2(b-a)}\int\limits_{a}^{b}%
\left| 2q(\frac{a+b}{2},s)+q(a,s)\right| ds.  \label{j14}
\end{equation}%
We now calculate%
\begin{eqnarray}
&&\int\limits_{a}^{b}\left| 2q(\frac{a+b}{2},s)+q(a,s)\right| ds  \label{j15}
\\
&=&\int\limits_{a}^{\frac{a+b}{2}}\left| -2\frac{b-a}{2}(s-a)^{2}-\frac{%
(b-a)^{2}}{2}(s-b)-\frac{b-a}{2}(s-b)^{2}\right| ds  \notag \\
&&+\int\limits_{\frac{a+b}{2}}^{b}\left| -2\frac{b-a}{2}(s-b)^{2}-\frac{%
(b-a)^{2}}{2}(s-b)-\frac{b-a}{2}(s-b)^{2}\right| ds.  \notag
\end{eqnarray}%
From the equation%
\begin{equation*}
-2\frac{b-a}{2}(s-a)^{2}-\frac{(b-a)^{2}}{2}(s-b)-\frac{b-a}{2}(s-b)^{2}=0
\end{equation*}%
we find the solutions%
\begin{equation}
s_{1}=a,\quad s_{2}=\frac{2a+b}{3}.  \label{j16}
\end{equation}%
From the equation%
\begin{equation*}
-2\frac{b-a}{2}(s-b)^{2}-\frac{(b-a)^{2}}{2}(s-b)-\frac{b-a}{2}(s-b)^{2}=0
\end{equation*}%
we find the solutions%
\begin{equation}
s_{3}=b,\quad s_{4}=\frac{a+2b}{3}.  \label{j17}
\end{equation}%
From (\ref{j15})-(\ref{j17}) we have%
\begin{eqnarray*}
&&\int\limits_{a}^{b}\left| 2q(\frac{a+b}{2},s)+q(a,s)\right| ds \\
&=&\int\limits_{a}^{\frac{2a+b}{3}}\left[ -2\frac{b-a}{2}(s-a)^{2}-\frac{%
(b-a)^{2}}{2}(s-b)-\frac{b-a}{2}(s-b)^{2}\right] ds \\
&&+\int\limits_{\frac{2a+b}{3}}^{\frac{a+b}{2}}\left[ 2\frac{b-a}{2}%
(s-a)^{2}+\frac{(b-a)^{2}}{2}(s-b)+\frac{b-a}{2}(s-b)^{2}\right] ds \\
&&+\int\limits_{\frac{a+b}{2}}^{\frac{a+2b}{3}}\left[ 2\frac{b-a}{2}%
(s-b)^{2}+\frac{(b-a)^{2}}{2}(s-b)+\frac{b-a}{2}(s-b)^{2}\right] ds \\
&&+\int\limits_{\frac{a+2b}{3}}^{b}\left[ -2\frac{b-a}{2}(s-b)^{2}-\frac{%
(b-a)^{2}}{2}(s-b)-\frac{b-a}{2}(s-b)^{2}\right] ds \\
&=&4\frac{1}{108}(b-a)^{4}=\frac{1}{27}(b-a)^{4}.
\end{eqnarray*}%
From the above relation and (\ref{j11a})--(\ref{j14}) we get%
\begin{equation}
\left| \frac{b-a}{2}\left[ f(a)+4f(\frac{a+b}{2})+f(b)\right]
-3\int\limits_{a}^{b}f(t)dt\right| \leq \frac{\Gamma -\gamma }{54}(b-a)^{3}.
\label{j19}
\end{equation}%
From (\ref{j19}) we easily get (\ref{j7}).
\end{proof}

\begin{remark}
\label{R1}The usual Peano error bound is%
\begin{equation}
\left| \frac{b-a}{6}\left[ f(a)+4f(\frac{a+b}{2})+f(b)\right]
-\int\limits_{a}^{b}f(t)dt\right| \leq \frac{(b-a)^{3}}{81}\left\| f^{\prime
\prime }\right\| _{\infty }.  \label{p1}
\end{equation}%
If we choose%
\begin{equation*}
\gamma =\underset{s\in \left[ a,b\right] }{\inf }f^{\prime \prime }(s)\text{%
, \ }\Gamma =\underset{s\in \left[ a,b\right] }{\sup }f^{\prime \prime }(s)
\end{equation*}%
then $\frac{\Gamma -\gamma }{2}\leq \left\| f^{\prime \prime }\right\|
_{\infty }$ and it is obvious that (\ref{j7}) is better than (\ref{p1}). In
fact, these two bounds are equal if and only if $\Gamma =-\gamma .$ This
case ($\Gamma =-\gamma $) is very rare in practice. Specially, if $\Gamma $
is large and $\Gamma \approx \gamma $ then (\ref{j7}) is much better than (%
\ref{p1}).
\end{remark}

\section{3/8 Simpson rule}

\begin{theorem}
\label{T2}Under the assumptions of Theorem \ref{T1} we have%
\begin{equation}
\left| \frac{b-a}{8}\left[ f(a)+3f(\frac{2a+b}{3})+3f(\frac{a+2b}{3})+f(b)%
\right] -\int\limits_{a}^{b}f(t)dt\right| \leq \frac{\Gamma -\gamma }{384}%
(b-a)^{3}.  \label{j21}
\end{equation}
\end{theorem}

\begin{proof}
If we substitute $x=\frac{2a+b}{3}$ in (\ref{j1}) then we get the
corresponding left-hand side%
\begin{equation}
\frac{b-a}{6}\left[ f(b)-f(a)\right] +f(\frac{2a+b}{3})(b-a)-\int%
\limits_{a}^{b}f(t)dt.  \label{j22}
\end{equation}%
For $x=\frac{a+2b}{3}$ we have the corresponding left-hand side%
\begin{equation}
-\frac{b-a}{6}\left[ f(b)-f(a)\right] +f(\frac{a+2b}{3})(b-a)-\int%
\limits_{a}^{b}f(t)dt.  \label{j23}
\end{equation}%
If we now multiply (\ref{j22}) and (\ref{j23}) by 3 and add (\ref{j8}) and (%
\ref{j9}) then we get the left-hand side of the form%
\begin{equation}
(b-a)\left[ f(a)+3f(\frac{2a+b}{3})+3f(\frac{a+2b}{3})+f(b)\right]
-8\int\limits_{a}^{b}f(t)dt.  \label{j24}
\end{equation}%
The corresponding right-hand side is%
\begin{equation}
R(a,b)=\frac{1}{b-a}\int\limits_{a}^{b}\left[ q(a,s)+3q(\frac{2a+b}{3},s)+3q(%
\frac{a+2b}{3},s)+q(b,s)\right] f^{\prime \prime }(s)ds.  \label{j25}
\end{equation}%
Let $K_{2}(s)=q(a,s)+3q(\frac{2a+b}{3},s)+3q(\frac{a+2b}{3},s)+q(b,s)$. Then
we have%
\begin{equation*}
\int\limits_{a}^{b}K_{2}(s)ds=0
\end{equation*}%
such that%
\begin{equation*}
\frac{1}{b-a}\int\limits_{a}^{b}K_{2}(s)f^{\prime \prime }(s)ds=\frac{1}{b-a}%
\int\limits_{a}^{b}K_{2}(s)\left[ f^{\prime \prime }(s)-\frac{\Gamma +\gamma 
}{2}\right] ds
\end{equation*}%
and%
\begin{eqnarray}
\left| R(a,b)\right| &\leq &\frac{1}{b-a}\underset{s\in \left[ a,b\right] }{%
\max }\left| f^{\prime \prime }(s)-\frac{\Gamma +\gamma }{2}\right|
\int\limits_{a}^{b}\left| K_{2}(s)\right| ds  \label{ub2} \\
&\leq &\frac{\Gamma -\gamma }{2(b-a)}\int\limits_{a}^{b}\left|
K_{2}(s)\right| ds,  \notag
\end{eqnarray}%
since (\ref{zvj}) holds. We now calculate%
\begin{eqnarray*}
&&\int\limits_{a}^{b}\left| q(a,s)+3q(\frac{2a+b}{3},s)+3q(\frac{a+2b}{3}%
,s)+q(b,s)\right| ds \\
&=&\int\limits_{a}^{\frac{2a+b}{3}}\left| \frac{1}{2}(b-a)^{3}-\frac{1}{2}%
(b-a)(s-b)^{2}-\frac{7}{2}(b-a)(s-a)^{2}\right| ds \\
&&+\int\limits_{\frac{2a+b}{3}}^{\frac{a+2b}{3}}\left| 2(b-a)(\left[
(s-a)^{2}+(s-b)^{2}\right] -(b-a)^{3}\right| ds \\
&&+\int\limits_{\frac{a+2b}{3}}^{b}\left| \frac{1}{2}(b-a)^{3}-\frac{1}{2}%
(b-a)(s-a)^{2}-\frac{7}{2}(b-a)(s-b)^{2}\right| ds.
\end{eqnarray*}%
From the equation%
\begin{equation*}
\frac{1}{2}(b-a)^{3}-\frac{1}{2}(b-a)(s-b)^{2}-\frac{7}{2}(b-a)(s-a)^{2}=0.
\end{equation*}%
we find the solutions%
\begin{equation}
s_{1}=a,\quad s_{2}=\frac{3a+b}{4}.  \label{j210}
\end{equation}%
From the equation%
\begin{equation*}
2(b-a)(\left[ (s-a)^{2}+(s-b)^{2}\right] -(b-a)^{3}=0
\end{equation*}%
we find the solution%
\begin{equation}
s_{3}=\frac{a+b}{2}.  \label{j211}
\end{equation}%
From the equation%
\begin{equation*}
\frac{1}{2}(b-a)^{3}-\frac{1}{2}(b-a)(s-a)^{2}-\frac{7}{2}(b-a)(s-b)^{2}=0
\end{equation*}%
we find the solutions%
\begin{equation}
s_{4}=b,\quad s_{5}=\frac{a+3b}{4}.  \label{j212}
\end{equation}%
From the above relations we get%
\begin{eqnarray*}
&&\int\limits_{a}^{b}\left| q(a,s)+3q(\frac{2a+b}{3},s)+3q(\frac{a+2b}{3}%
,s)+q(b,s)\right| ds \\
&=&\int\limits_{a}^{\frac{3a+b}{4}}\left[ \frac{1}{2}(b-a)^{3}-\frac{1}{2}%
(b-a)(s-b)^{2}-\frac{7}{2}(b-a)(s-a)^{2}\right] ds \\
&&-\int\limits_{\frac{3a+b}{4}}^{\frac{2a+b}{3}}\left[ \frac{1}{2}(b-a)^{3}-%
\frac{1}{2}(b-a)(s-b)^{2}-\frac{7}{2}(b-a)(s-a)^{2}\right] ds \\
&&+\int\limits_{\frac{2a+b}{3}}^{\frac{a+2b}{3}}\left[
2(b-a)((s-b)^{2}+(s-a)^{2})-(b-a)^{3}\right] ds \\
&&-\int\limits_{\frac{a+2b}{3}}^{\frac{a+3b}{4}}\left[ \frac{1}{2}(b-a)^{3}-%
\frac{1}{2}(b-a)(s-a)^{2}-\frac{7}{2}(b-a)(s-b)^{2}\right] ds \\
&&+\int\limits_{\frac{a+3b}{4}}^{b}\left[ \frac{1}{2}(b-a)^{3}-\frac{1}{2}%
(b-a)(s-a)^{2}-\frac{7}{2}(b-a)(s-b)^{2}\right] ds \\
&=&\frac{1}{24}(b-a)^{4}.
\end{eqnarray*}%
From the above relation and (\ref{j24}) and (\ref{ub2}) we get%
\begin{equation*}
\left| (b-a)\left[ f(a)+3f(\frac{2a+b}{3})+3f(\frac{a+2b}{3})+f(b)\right]
-8\int\limits_{a}^{b}f(t)dt\right| \leq \frac{\Gamma -\gamma }{48}(b-a)^{3}.
\end{equation*}%
This completes the proof.
\end{proof}

\begin{remark}
\label{R2} The usual Peano error bound is%
\begin{equation*}
\left| \frac{b-a}{8}\left[ f(a)+3f(\frac{2a+b}{3})+3f(\frac{a+2b}{3})+f(b)%
\right] -\int\limits_{a}^{b}f(t)dt\right| \leq \frac{\left\| f^{\prime
\prime }\right\| _{\infty }}{192}(b-a)^{3}.
\end{equation*}%
For the reasons given in Remark \ref{R1} the estimation obtained in Theorem %
\ref{T2}, which is a Peano-like bound, is better than the above Peano bound.
\end{remark}

\section{Boole's rule}

\begin{theorem}
\label{T3}Under the assumptions of Theorem \ref{T2} we have%
\begin{equation}
\begin{array}{c}
\left| \frac{b-a}{90}\left[ 7f(a)+32f(\frac{3a+b}{4})+12f(\frac{a+b}{2})+32f(%
\frac{a+3b}{4})+7f(b)\right] \right. \\ 
\left. -\int\limits_{a}^{b}f(t)dt\right| \leq \frac{509}{273\,375}(\Gamma
-\gamma )(b-a)^{3}.%
\end{array}
\label{j320}
\end{equation}
\end{theorem}

\begin{proof}
We first write left-hand sides of (\ref{j1}) for $x=\frac{3a+b}{4}$ and $x=%
\frac{a+3b}{4}.$ We have%
\begin{equation}
f(\frac{3a+b}{4})+\frac{b-a}{4}\left[ f(b)-f(a)\right] -\int%
\limits_{a}^{b}f(t)dt  \label{j321}
\end{equation}%
and%
\begin{equation}
f(\frac{a+3b}{4})-\frac{b-a}{4}\left[ f(b)-f(a)\right] -\int%
\limits_{a}^{b}f(t)dt.  \label{j322}
\end{equation}%
If we now multiply (\ref{j8}) and (\ref{j9}) by 7, (\ref{j10}) by 12, (\ref%
{j321}) and (\ref{j322}) by 32 and sum the obtained results then we get the
left-hand side of the form%
\begin{equation}
(b-a)\left[ 7f(a)+32f(\frac{3a+b}{4})+12f(\frac{a+b}{2})+32f(\frac{a+3b}{4}%
)+7f(b)\right] -90\int\limits_{a}^{b}f(t)dt.  \label{j323}
\end{equation}%
For the corresponding right-hand side we get%
\begin{eqnarray}
&&R(a,b)  \label{j327} \\
&=&\frac{1}{b-a}\int\limits_{a}^{b}\left[ 7q(a,s)+32q(\frac{3a+b}{4},s)+12q(%
\frac{a+b}{2},s)+32q(\frac{a+3b}{4},s)\right.  \notag \\
&&\left. +7q(b,s)\right] f^{\prime \prime }(s)ds.  \notag
\end{eqnarray}

Let $K_{3}(s)=7q(a,s)+32q(\frac{3a+b}{4},s)+12q(\frac{a+b}{2},s)+32q(\frac{%
a+3b}{4},s)+7q(b,s)$. Then we have%
\begin{equation*}
\int\limits_{a}^{b}K_{3}(s)ds=0
\end{equation*}%
such that%
\begin{equation*}
\frac{1}{b-a}\int\limits_{a}^{b}K_{3}(s)f^{\prime \prime }(s)ds=\frac{1}{b-a}%
\int\limits_{a}^{b}K_{3}(s)\left[ f^{\prime \prime }(s)-\frac{\Gamma +\gamma 
}{2}\right] ds
\end{equation*}%
and%
\begin{eqnarray*}
\left| R(a,b)\right| &\leq &\frac{1}{b-a}\underset{s\in \left[ a,b\right] }{%
\max }\left| f^{\prime \prime }(s)-\frac{\Gamma +\gamma }{2}\right|
\int\limits_{a}^{b}\left| K_{3}(s)\right| ds \\
&\leq &\frac{\Gamma -\gamma }{2(b-a)}\int\limits_{a}^{b}\left|
K_{3}(s)\right| ds,
\end{eqnarray*}%
since (\ref{zvj}) holds. Hence,%
\begin{eqnarray*}
\left| R(a,b)\right| &\leq &\frac{\Gamma -\gamma }{2(b-a)}%
\int\limits_{a}^{b}\left| 7q(a,s)+32q(\frac{3a+b}{4},s)+12q(\frac{a+b}{2}%
,s)+32q(\frac{a+3b}{4},s)\right. \\
&&\left. +7q(b,s)\right| ds.
\end{eqnarray*}%
We now calculate%
\begin{eqnarray*}
&&\int\limits_{a}^{b}\left| 7q(a,s)+32q(\frac{3a+b}{4},s)+12q(\frac{a+b}{2}%
,s)+32q(\frac{a+3b}{4},s)+7q(b,s)\right| ds \\
&=&\int\limits_{a}^{\frac{3a+b}{4}}\left| -\frac{7(b-a)}{2}\left[
(s-b)^{2}+(s-a)^{2}\right] -38(b-a)(s-a)^{2}+\frac{7}{2}(b-a)^{3}\right| ds
\\
&&+\int\limits_{\frac{3a+b}{4}}^{\frac{a+b}{2}}\left| 39(b-a)\left[ \frac{%
(b-a)^{2}}{4}+(s-\frac{a+b}{2})^{2}\right] -\frac{23}{2}%
(b-a)^{3}+6(b-a)(s-a)^{2}\right| ds \\
&&+\int\limits_{\frac{a+b}{2}}^{\frac{a+3b}{4}}\left| 39(b-a)\left[ \frac{%
(b-a)^{2}}{4}+(s-\frac{a+b}{2})^{2}\right] -\frac{23}{2}%
(b-a)^{3}+6(b-a)(s-b)^{2}\right| ds \\
&&+\int\limits_{\frac{a+3b}{4}}^{b}\left| -\frac{7(b-a)}{2}\left[
(s-b)^{2}+(s-a)^{2}\right] -38(b-a)(s-b)^{2}+\frac{7}{2}(b-a)^{3}\right| ds.
\end{eqnarray*}%
Let us denote the integrands of the right-hand side of the above relation by 
$Q_{1}(s)$, $Q_{2}(s)$, $Q_{3}(s)$ and $Q_{4}(s)$, respectively. These
integrands have the next zero points%
\begin{equation*}
s_{0}=a,\text{ }s_{1}=\frac{38a+7b}{45},\text{ }s_{2}=\frac{2a+b}{3},\text{ }%
s_{3}=\frac{a+2b}{3},\text{ }s_{4}=\frac{7a+38b}{45},\text{ }s_{5}=b.
\end{equation*}%
From the above two relations we get%
\begin{eqnarray*}
&&\int\limits_{a}^{b}\left| 7q(a,s)+32q(\frac{3a+b}{4},s)+12q(\frac{a+b}{2}%
,s)+32q(\frac{a+3b}{4},s)+7q(b,s)\right| ds \\
&=&\int\limits_{a}^{s_{1}}Q_{1}(s)ds-\int\limits_{s_{1}}^{\frac{3a+b}{4}%
}Q_{1}(s)ds-\int\limits_{\frac{3a+b}{4}}^{s_{2}}Q_{2}(s)ds+\int%
\limits_{s_{2}}^{\frac{a+b}{2}}Q_{2}(s)ds \\
&&+\int\limits_{\frac{a+b}{2}}^{s_{3}}Q_{3}(s)ds-\int\limits_{s_{3}}^{\frac{%
a+3b}{4}}Q_{3}(s)ds-\int\limits_{\frac{a+3b}{4}}^{s_{4}}Q_{4}(s)ds+\int%
\limits_{s_{4}}^{b}Q_{4}(s)ds \\
&=&\frac{2036}{6075}(b-a)^{4}.
\end{eqnarray*}%
Thus, we have%
\begin{equation}
\left| R(a,b)\right| \leq \frac{1018}{6075}(\Gamma -\gamma )(b-a)^{3}.
\label{j331}
\end{equation}%
From (\ref{j323}), (\ref{j327}) and (\ref{j331}) we easily get (\ref{j320}).
\end{proof}

\begin{remark}
\label{R3} The usual Peano error bound is%
\begin{equation*}
\begin{array}{c}
\left| \frac{b-a}{90}\left[ 7f(a)+32f(\frac{3a+b}{4})+12f(\frac{a+b}{2})+32f(%
\frac{a+3b}{4})+7f(b)\right] \right. \\ 
\left. -\int\limits_{a}^{b}f(t)dt\right| \leq \frac{1018}{273\,375}\left\|
f^{\prime \prime }\right\| _{\infty }(b-a)^{3}.%
\end{array}%
\end{equation*}%
For the reasons given in Remark \ref{R1} the estimation obtained in Theorem %
\ref{T3}, which is a Peano-like bound, is better than the above Peano bound.
\end{remark}

\end{document}